\documentclass[12pt,reqno]{amsart}
\author{Kai (Steve) Fan}
\author{Paul Pollack}
\address{Department of Mathematics\\ University of Georgia\\ Athens, GA 30602}
\email{Steve.Fan@uga.edu}
\email{pollack@uga.edu}
\subjclass[2020]{}
\usepackage{mathtools}
\usepackage{empheq}
\usepackage{colonequals}
\usepackage[symbol]{footmisc}

\usepackage{amsmath,amssymb,amsthm,geometry,mathscinet,enumitem,url,amsrefs}
\usepackage{parskip}
\usepackage{stackengine}
\usepackage{colonequals}
\usepackage{xcolor}

\renewcommand\supset\supseteq
\renewcommand\phi\varphi
\usepackage[utf8]{inputenc}
\makeatletter
\NewCommandCopy\@@pmod\pmod
\DeclareRobustCommand{\pmod}{\@ifstar\@pmods\@@pmod}
\def\@pmods#1{\mkern4mu({\operator@font mod}\mkern 6mu#1)}
\makeatother

\numberwithin{equation}{section}

\DeclareMathAlphabet{\curly}{U}{rsfs}{m}{n}
\newtheorem{thm}{Theorem}[section]

\newtheorem{prop}[thm]{Proposition}

\theoremstyle{remark}

\usepackage{enumitem}

\def\N{\mathbb{N}}

\def\Dd{\mathcal{D}}

\def\Z{\mathbb{Z}}

\def\Dd{\mathcal{D}}

\renewcommand\subset\subseteq
\newcommand\Li{\mathrm{Li}}

\begin{document}

\title{The maximal order of the shifted-prime divisor function}
\dedicatory{Dedicated to the 80\textsuperscript{\,th} Birthdays of 
Melvyn Nathanson and Carl Pomerance}
\begin{abstract} 
For each positive integer $n$, we denote by $\omega^*(n)$ the number of
shifted-prime divisors $p-1$ of $n$, i.e., 
\[\omega^*(n)\colonequals\sum_{p-1\mid n}1.\]
First introduced by Prachar in 1955, this function has interesting applications in primality testing and bears a strong connection with counting Carmichael numbers. Prachar showed that for a certain constant $c_0 > 0$,
\[\omega^*(n)>\exp\left(c_0\frac{\log n}{(\log\log n)^2}\right)\]
for infinitely many $n$. This result was later improved by Adleman, Pomerance and Rumely, who established an inequality of the same shape  with $(\log\log n)^2$ replaced by $\log\log n$. Assuming the Generalized Riemann Hypothesis for Dirichlet $L$-functions, Prachar also proved the stronger inequality 
\[\omega^*(n)>\exp\left(\left(\frac{1}{2}\log2+o(1)\right)\frac{\log n}{\log\log n}\right)\]
for infinitely many $n$. By refining the arguments of Prachar and of Adleman, Pomerance and Rumely, we improve on their results by establishing 
\begin{align*}
\omega^*(n)&>\exp\left(0.6736\log 2\cdot\frac{\log n}{\log\log n}\right) \text{~~~(unconditionally)},\\
\omega^*(n)&>\exp\left(\left(\log\left(\frac{1+\sqrt{5}}{2}\right)+o(1)\right)\frac{\log n}{\log\log n}\right) \text{~~~(under GRH)},
\end{align*}
for infinitely many $n$.
\end{abstract}

\maketitle

\section{Introduction}
For each $n\in\N$, let $\omega^*(n)$ count the number of
shifted primes of the form $p-1$ dividing $n$, that is, 
\[\omega^*(n)\colonequals\sum_{p-1\mid n}1.\]
This function has found interesting applications in primality testing \cite{APR83} and is closely related to counting Carmichael numbers \cite{AGP94}. It is easy to see that $\omega^*(n)\ge1$ with equality if and only if $n$ is odd. On the other hand, the function $\omega^*$ can attain fairly large values infinitely often. Indeed, Prachar, who initiated the study of $\omega^*$ in his influential work \cite{Prachar55}, showed that there is some absolute constant $c_1>0$ such that
\begin{equation}\label{eq:Prachar}
\omega^*(n)>\exp\left(c_1\frac{\log n}{(\log\log n)^2}\right)
\end{equation}
for infinitely many $n$. Assuming the Generalized Riemann Hypothesis for Dirichlet $L$-functions (GRH), he was able to save a $\log\log n$ factor and provide a lower bound for $c_1$. More precisely, he proved that GRH implies the inequality
\begin{equation}\label{eq:PracharGRH}
\omega^*(n)>\exp\left(\left(\frac{1}{2}\log2+o(1)\right)\frac{\log n}{\log\log n}\right)
\end{equation}
for infinitely many $n$. Nearly three decades later, Adleman, Pomerance and Rumely \cite{APR83} improved Prachar's unconditional lower bound \eqref{eq:Prachar} to
\begin{equation}\label{eq:APR}
\omega^*(n)>\exp\left(c_2\frac{\log n}{\log\log n}\right),
\end{equation}
along a sequence of $n$ tending to infinity, where $c_2>0$ is some absolute constant. This lower bound has the same shape as Prachar's GRH-conditional lower bound \eqref{eq:PracharGRH}, up to the constant factor. They also conjectured that the constant $\frac{1}{2}\log 2$ in \eqref{eq:PracharGRH} can in fact be upgraded to $\log2$. It is worth noting that the conjectured constant $\log 2$ would be best possible, since 
\[\omega^*(n)\le \tau(n)\le \exp\left((\log2+o(1))\frac{\log n}{\log\log n}\right)\]
for sufficiently large $n$ \cite[Theorem 317]{HW08}, where $\tau(n)$ denotes the number of positive divisors of $n$. Interestingly, this conjecture, if true, would imply that the Carmichael function $\lambda(n)$, which may be defined as the exponent of the multiplicative group $(\Z/n\Z)^{\times}$, satisfies 
\[\liminf_{n\to\infty}\frac{\log\lambda(n)}{(\log\log n)\log\log\log n}=\frac{1}{\log 2},\]
an observation first made by Erd\H{o}s, Pomerance, and Schmutz \cite{EPS91}.

The study of $\omega^*$ was recently revived by Murty and Murty \cite{MM21}. For each $k\in\N$, we define the $k$th moment of $\omega^*$ by
\[M_k(x)\colonequals\frac{1}{x}\sum_{n\le x}\omega^*(n)^k.\]
These moments encode useful information about the distribution of $\omega^*$. A quick application of Mertens' theorem \cite[Theorem 427]{HW08} yields $M_1(x)=\log\log x+O(1)$. In the same paper \cite{Prachar55}, Prachar showed that $M_2(x)\ll(\log x)^2$. In their recent work \cite{MM21}, Murty and Murty \cite{MM21} proved the estimate $(\log\log x)^3\ll M_2(x)\ll\log x$ and conjectured the asymptotic formula $M_2(x)\sim C_2\log x$ with some constant $C_2>0$. Two years later, Ding \cite{Ding23} achieved the order matching lower bound $M_2(x)\gg \log x$. More recently, Pomerance and the first author \cite{FP24} established the estimate $M_3(x)\asymp(\log x)^4$ and made the more general conjecture that for each $k\ge2$ there exists a constant $C_k>0$ such that $M_k(x)\sim C_k(\log x)^{2^k-k-1}$. Based on a heuristic argument, the first author \cite{Fan25} also conjectured that $C_2=\zeta(2)^2\zeta(3)/\zeta(6)$, where $\zeta$ is the Riemann zeta-function. Despite the fact that no asymptotic formula is currently known even for a single value of $k\ge2$, quite recently Gabdullin \cite{Gab25} pinned down the correct order of magnitude for $M_k(x)$, showing that $M_k(x)\asymp(\log x)^{2^k-k-1}$ for every fixed $k\ge2$.

The main objective of this paper is to furnish explicit values for the constant $c_2$ appearing in \eqref{eq:APR}. Perhaps a bit surprisingly, we are able to obtain unconditionally a numerical value for $c_2$ which is slightly larger than $\frac{2}{3}\log 2$, hence surpassing the constant $\frac{1}{2}\log 2$ in Prachar's GRH-conditional lower bound \eqref{eq:PracharGRH}. Under GRH, our argument yields a numerical value for $c_2$ which exceeds $(\log 2)^2$. The following theorem provides a precise statement of our results.

\begin{thm}\label{thm:maxomega*}
There exist infinitely many $n$ such that
\begin{equation}\label{eq:maxomega*}
\omega^*(n)>\exp\left(0.6736\log 2\cdot\frac{\log n}{\log\log n}\right).    
\end{equation}
Moreover, if GRH is true, then we have
\begin{equation}\label{eq:GRHmaxomega*}
\omega^*(n)>\exp\left(\left(\log\left(\frac{1+\sqrt{5}}{2}\right)+o(1)\right)\frac{\log n}{\log\log n}\right)   
\end{equation}
for infinitely many $n$.
\end{thm}

% \textcolor{blue}{We got 0.6269 unconditionally using Thm 2.1 in AGP94, but maybe we can improve this by using the analogous result of Harman which allows larger moduli in progressions but with a slightly more complicated exceptional set.}

% \textcolor{blue}{Do we also want to include a discussion on the connection of this problem with known conjectures such as Montgomery's conjecture on primes in progressions and Carl's conjecture on smooth shifted primes?}

% \textcolor{magenta}{Yes, I think it would be nice to record some version of the arguments we discussed by email, at the end of the paper. That Montgomery's conjecture gives the correct constant will be clear by that point. But the connection with conjectures on smooth (or powersmooth, or squarefree) shifted primes seems worth a discussion.}

% \textcolor{blue}{Sounds good to me.}

\medskip
\section{An overview of Prachar's argument}\label{sec:overview}
We first outline Prachar's simple proof of \eqref{eq:PracharGRH} upon which our proof of Theorem \ref{thm:maxomega*} is built. In what follows, the letter $p$ always denotes a prime. We write $\pi(x)$ for the number of primes $p\le x$, and let $\pi(x;d,a)$ count the number of primes $p\le x$ with $p\equiv a\pmod*{d}$.

Suppose that $x$ is sufficiently large. Let $\epsilon\in(0,1)$ be arbitrary, and put
\begin{equation}\label{eq:choice of k}
k=\prod_{p\le (1-\epsilon)\log x}p=x^{1-\epsilon+o(1)}.\footnote[2]{In fact, Prachar \cite{Prachar55} took 
\[k=\prod_{p\le (1/2-\epsilon)\log x}p.\] 
However, this appears to be a misstep: This choice of $k$ only yields the smaller constant $\frac{1}{4}\log 2$ rather than $\frac{1}{2}\log 2$ asserted in \cite{Prachar55} and \eqref{eq:PracharGRH}.}   
\end{equation}
Under GRH we have \cite[Corollary 13.8]{MV06}
\[\pi(x;d,1)=\frac{\Li(x)}{\varphi(d)}+O\left(\sqrt{x}\log x\right)\gg \frac{x}{\varphi(d)\log x}\]
uniformly for all positive integers $d\le\sqrt{k}$, where
\[\Li(x)\colonequals\int_{2}^{x}\frac{\text{d}t}{\log t}.\]
For each $1\le d\le\sqrt{k}$ dividing $k$, denote by $A_d$ the number of pairs $(m,p)$, where $m\in\N$ and $p$ is prime, such that $m,p\le x$, 
\begin{equation}\label{eq:(m,p)}
 p\equiv 1\pmod*{d}\quad\text{and}\quad\gcd(m,k)=\frac{k}{d}.
\end{equation}
% \textcolor{magenta}{Maybe we just want to include the first two conditions in this equation? It seems the third one, with mod $k$ instead of mod $d$, is a consequence of the first two?} \textcolor{blue}{Agree. I have deleted the last one}. 
Let $A$ record the total number of pairs $(m,p)$ with $m,p\le x$ satisfying the congruence $m(p-1)\equiv 0\pmod*{k}$. Every pair counted by $A_d$, for some $d \mid k, d\le \sqrt{k}$, is counted by $A$. Moreover, an individual pair is counted by $A_d$ for at most one $d$. Since the number of choices for $p$ is clearly $\pi(x;d,1)$, and since the count of $m$ is at least $\lfloor x/k\rfloor\varphi(d)$, we have
\[A_d\ge\pi(x;d,1)\left\lfloor\frac{x}{k}\right\rfloor\varphi(d)\gg\frac{x^2}{k\log x},\]
from which it follows that
\[A\ge\sum_{\substack{d\le\sqrt{k}\\d\mid k}}A_d\gg\frac{x^2}{k\log x}\sum_{\substack{d\le\sqrt{k}\\d\mid k}}1 \ge\frac{\tau(k)x^2}{2k\log x}.\]
On the other hand, each pair $(m,p)$ counted by $A$ yields a positive integer $n=m(p-1)\le x^2$ divisible by $k$. Thus,
\[A\le\sum_{\substack{n\le x^2\\k\mid n}}\#\{(m,p)\colon m,p\le x\text{~and~}n=m(p-1)\}.\]
% {\color{magenta} I think we should probably describe how to easily see this. I can do this later if you don't get to it first.} \textcolor{blue}{I have modified the writing to make the argument clearer.}
Combining the above upper and lower bounds for $A$, we deduce that there exists some large $n\le x^2$ which is divisible by $k$ and admits
\[\gg\frac{\tau(k)x^2/(k\log x)}{x^2/k}=\frac{\tau(k)}{\log x}=\frac{2^{\pi((1-\epsilon)\log x)}}{\log x}>\exp\left(\left(\log2-2\epsilon\right)\frac{\log x}{\log\log x}\right)\]
representations of the form $n=m(p-1)$. For this particular $n$, we have
\[\omega^*(n)\ge\exp\left(\left(\log2-2\epsilon\right)\frac{\log x}{\log\log x}\right)>\exp\left(\left(\frac{1}{2}\log2-2\epsilon\right)\frac{\log n}{\log\log n}\right),\]
as desired.

In the argument outlined above, the representations $n=m(p-1)$ counted by $A_d$ all have $m,p\le x$ and $d\le\sqrt{k}$, with $k$ defined by \eqref{eq:choice of k}. In the next section, we shall prove Theorem \ref{thm:maxomega*} by adjusting the choices for $k,d,m,p$. We use insights from the theory of anatomy of integers to maximize the total number of representations $n=m(p-1)$ counted by $A$.

\medskip
\section{Refining Prachar's argument: Proof of Theorem \ref{thm:maxomega*}}
We start by proving the GRH-conditional inequality \eqref{eq:GRHmaxomega*}. Suppose that $x$ is sufficiently large. Put
\begin{equation}\label{eq:epsilonchoice} \epsilon \coloneqq (\log\log{x})^{-1/2}\quad\text{and}\quad u \coloneqq \frac{3+\sqrt{5}}{4},  \end{equation}
%Then $\sqrt{2u}=\big(1+\sqrt{5}\big)/2$ coincides with the golden ratio.
%Define
and set
\begin{align*} k&\coloneqq \prod_{p\le (u-\epsilon)\log x}p \\
&=x^{u-\epsilon}\exp\left(O\left(\frac{\log{x}}{\log\log{x}}\right)\right).\end{align*} 

We will show momentarily that there is a set $\mathcal{D}$ of divisors of $k$, where each $d \in \mathcal{D}$ obeys the estimate
\begin{equation}\label{eq:destimate} d= x^{\frac{1}{2}-\epsilon} \exp\left(O\left(\frac{\log{x}}{(\log\log{x})^2}\right)\right), \end{equation}
and where, with $X\coloneqq x^{\frac{1}{2}+u}$, 
\begin{equation}\label{eq:sizedlower} \#\Dd \ge \exp\left(\left(\log\left(\frac{1+\sqrt{5}}{2}\right)+o(1)\right)\frac{\log{X}}{\log\log{X}}\right). \end{equation}
Let us see now how this claim implies \eqref{eq:GRHmaxomega*}.

For each $d\in \Dd$, let $A_d$ denote the number of pairs $(m,p)$ with $m\le y_d\colonequals x^{u}/d$ and $p\le x$, satisfying the same conditions described in \eqref{eq:(m,p)}. By inclusion-exclusion, the number of choices for $m$ is 
\[\sum_{\substack{m'\le y_d/(k/d)\\(m',d)=1}}1=\frac{\varphi(d)}{d}\cdot\frac{y_d}{k/d}+O(\tau(d))\gg\frac{\varphi(d)}{d}\cdot\frac{x^u}{k},\]
since $\tau(d)\le \exp\left(O\left(\frac{\log{x}}{\log\log{x}}\right)\right)$, and 
\[\frac{\varphi(d)}{d}\cdot\frac{y_d}{k/d}=\frac{\varphi(d)}{d}\cdot\frac{x^u}{k}\gg \exp\left(\frac{\log{x}}{2(\log\log{x})^{1/2}}\right).\]
Our estimate \eqref{eq:destimate} for $d$, along with the choice of $\epsilon$ in \eqref{eq:epsilonchoice}, implies that $d \le x^{1/2}/(\log{x})^3$ (say), so that 
\[\pi(x;d,1)=\frac{\Li(x)}{\varphi(d)}+O\left(\sqrt{x}\log x\right)\gg \frac{x}{\varphi(d)\log x}. \]

It follows that for all $d \in \Dd$, 
\begin{equation}\label{eq:A_d}
A_d\gg\pi(x;d,1)\frac{\varphi(d)}{d}\cdot\frac{x^u}{k}\gg\frac{x^{u+1}}{kd\log x} = \frac{x^{u+\frac12 + \epsilon}}{k\log{x}} \exp\left(O\left(\frac{\log{x}}{(\log\log{x})^2}\right)\right).
\end{equation}
Therefore, 
\begin{equation}\label{eq:tocompare0} \sum_{d \in \mathcal{D}} A_d \gg \#\mathcal{D}\,\frac{x^{u+\frac12 + \epsilon}}{k\log{x}} \exp\left(O\left(\frac{\log{x}}{(\log\log{x})^2}\right)\right). \end{equation}
On the other hand, if $(m,p)$ is counted by some $A_d$, then $m(p-1)$ is a multiple of $k$ and \begin{equation}\label{eq:tocompare1} m(p-1) \le \big(\max_{d \in \mathcal{D}} y_d\big) x = x^{u + \frac{1}{2}+\epsilon} \exp\left(O\left(\frac{\log{x}}{(\log\log{x})^2}\right)\right).\end{equation}
Reasoning as in \S\ref{sec:overview}, we conclude upon comparing \eqref{eq:tocompare0} and \eqref{eq:tocompare1} that there is a large multiple $n$ of $k$, not exceeding the final expression in \eqref{eq:tocompare1}, with 
\[ \omega^{*}(n) \ge \frac{\#\mathcal{D}}{\log{x}} \exp\left(O\left(\frac{\log{x}}{(\log\log{x})^2}\right)\right). \]
Substituting the lower bound \eqref{eq:sizedlower} for $\#\Dd$ gives
\begin{align*} \omega^{*}(n) &\ge \exp\left(\left(\log\left(\frac{1+\sqrt{5}}{2}\right)+o(1)\right)\frac{\log{X}}{\log\log{X}}\right)\\&\ge \exp\left(\left(\log\left(\frac{1+\sqrt{5}}{2}\right)+o(1)\right)\frac{\log{n}}{\log\log{n}}\right), \end{align*}
where we use in the second line that $n \le X^{1+o(1)}$.

To show the existence of the set $\Dd$, we employ the probabilistic method. Let  \[ L = (u-\epsilon) \log{x}\quad\text{and} \quad R = \pi(L). \]
Furthermore,  set
\[ \rho = \frac{\frac{1}{2}-\epsilon}{u-\epsilon}. \]
We introduce i.i.d.\ Bernoulli random variables $v_r$, for each prime $r\le L$, where every $v_r$ takes the value $1$ with probability $\rho$. Then
\[ d \coloneqq \prod_{r \le L} r^{v_r} \]
is a random divisor of $k$. We proceed to study the distribution of $d$. 

It is straightforward to compute the expectation and variance of the random variable $\log{d}= \sum_{r \le L} v_r \log{r}$: We have
\[ \mathbb{E}[\log{d}] = \sum_{r \le L} \mathbb{E}[v_r \log{r}] = \rho \sum_{r \le L} \log{r} = \left(\frac{1}{2}-\epsilon\right)\log{x} + O(L/(\log{L})^3),\]
\[ \mathbb{V}[\log{d}] = \sum_{r \le L} \mathbb{V}[v_r \log{r}] = \rho(1-\rho) \sum_{r \le L} (\log{r})^2 \ll \sum_{r \le L} (\log{r})^2 \ll L\log{L}. \] This latter estimate implies, via Chebyshev's inequality, that $|\log{d} - \mathbb{E}[\log{d}]| > L^{2/3}$ with probability $O(L^{-1/3}\log{L}) = o(1)$. Thus, with probability $1+o(1)$,
\begin{align}\notag \left|\log{d}-\left(\frac{1}{2}-\epsilon\right)\log{x}\right| &\le \left|\log{d}-\mathbb{E}[\log{d}]\right| + \left|\mathbb{E}[\log{d}]- \left(\frac{1}{2}-\epsilon\right)\log{x}\right| 
\\ &\le L^{2/3} + L/(\log{L})^2 < 2L/(\log{L})^2.\label{eq:dclose}\end{align} 

Let $\mathcal{D}$ be the set of divisors $d$ of $k$ for which $|\log{d}-(\frac{1}{2}-\epsilon)\log{x}| < 2L/(\log{L})^2$. For each $d \in \mathcal{D}$, the desired estimate \eqref{eq:destimate} holds, and we have just seen that $\mathbb{P}(d \in \mathcal{D}) = 1+o(1)$, as $x\to\infty$. We proceed to translate this probability bound into a lower bound on $\#\Dd$. In fact, we will obtain the claimed lower bound \eqref{eq:sizedlower} for a certain convenient subset of $\mathcal{D}$.

The mean and variance of $\Omega(d)=\sum_{r \le L} v_r$ satisfy $\mathbb{E}[\Omega(d)] = \rho R$, $\mathbb{V}[\Omega(d)] \ll R$. 
% \textcolor{blue}{I have changed $\omega$ to $\Omega$ to make it agree with the standard usage of $\Omega(d)$ and avoid possible confusion with $\omega^*$.} 
So if we let $\mathcal{E}$ denote the set of $d\mid k$ with $|\Omega(d) - \rho R| > R^{2/3}$, then $\mathbb{P}(d \in \mathcal{E}) = o(1)$ by another application of Chebyshev's inequality. We put $\mathcal{D}' \coloneqq \mathcal{D}\setminus \mathcal{E}$ and observe that \[ \mathbb{P}(d \in \mathcal{D}') \ge \mathbb{P}(d \in \mathcal{D})- \mathbb{P}(d\in \mathcal{E}) = 1+o(1).\]

If $d \in \mathcal{D}'$, then $v_r=1$ for $\rho R+ O(R^{2/3})$ primes $r \le L$, while $v_r=0$ for $(1-\rho)R + O(R^{2/3})$ primes $r\le L$. \textbf{}
% \textcolor{blue}{I guess $X_r$ should be $v_r$?} 
Hence, each $d \in \mathcal{D}'$ carries a probability mass of
\[ \rho^{\rho R} (1-\rho)^{(1-\rho)R} \exp(O(R^{2/3})). \]
In order for the probability masses corresponding to $d\in\mathcal{D'}$ to sum to $1+o(1)$, it must be that
\begin{equation} \#\mathcal{D} \ge \#\mathcal{D'} \ge \rho^{-\rho R} (1-\rho)^{-(1-\rho)R} \exp(O(R^{2/3})). \label{eq:dboundmore0}\end{equation}
Since $R = (1+o(1)) u \frac{\log{x}}{\log\log{x}}$ while $\rho = \frac{1}{2u} + o(1)$, we have
\begin{equation}\label{eq:dboundmore1} \rho^{-\rho R} (1-\rho)^{-(1-\rho)R} \exp(O(R^{2/3})) = \exp((C+o(1))\log{x}/\log\log{x}), \end{equation}
where
\begin{align*} C &= \frac{1}{2}\log{(2u)} + \left(u-\frac{1}{2}\right) \log\frac{2u}{2u-1} \\
&= \left(u+\frac{1}{2}\right) \log\frac{2u}{2u-1} + \left(\log\sqrt{2u} - \log \frac{2u}{2u-1}\right)\\
&= \left(u+\frac{1}{2}\right) \log\frac{1+\sqrt{5}}{2},\end{align*}
noting for the last line that $\frac{2u}{2u-1} = \frac{1+\sqrt{5}}{2} = \sqrt{2u}$. Finally, since $X=x^{\frac{1}{2}+u}$, \begin{equation}\label{eq:dboundmore2} \left(u+\frac{1}{2}\right)\frac{\log{x}}{\log\log{x}} = (1+o(1)) \frac{\log{X}}{\log\log{X}}. \end{equation} The lower bound \eqref{eq:sizedlower} on $\#\mathcal{D}$ follows from \eqref{eq:dboundmore0}, \eqref{eq:dboundmore1} and \eqref{eq:dboundmore2}.

The unconditional inequality \eqref{eq:maxomega*} follows in a similar fashion, using the following result of Harman (see \cite[Theorem 1.2]{harman08}) as a proxy for the GRH. 

\begin{prop}\label{prop:harman} There is an absolute constant $\delta > 0$ making the following true.
  
For each $\eta > 0$, there are constants $K \ge 2$ and $c > 0$ such that the following holds. Suppose  
\[ K < d < x^{0.4736}, \quad\text{and}\quad p\mid d\Rightarrow p < d^{\delta}. \] 
Furthermore, assume that for every $f\mid d$ and primitive character $\chi\bmod{f}$,
\begin{equation}\label{eq:harmanregion} L(s,\chi) \ne 0\quad\text{for}\quad \mathrm{Re}(s) > 1-\frac{1}{(\log{d})^{3/4}},~ |\mathrm{Im}(s)|\le \exp(\eta(\log{d})^{3/4}). \end{equation}
Then for every $a$ with $\gcd(a,d)=1$, we have
\[ \pi(x;d,a) \ge \frac{cx}{\phi(d)\log{x}}. \]
\end{prop}

Put 
% \[\epsilon \coloneqq (\log\log{x})^{-1/2}\quad\text{and}\quad\theta\coloneqq 0.416666<\frac{5}{12},\]
% and define $u\approx1.182249$ to be the unique point at which the function
% \begin{equation}\label{eq:f_{theta}}
% f_{\theta}(t)\colonequals\frac{1}{t+1-\theta}\left(\theta\log\frac{t}{\theta}-\left(t-\theta\right)\log\left(1-\frac{\theta}{t}\right)\right)=\frac{t\log t-(t-\theta)\log(t-\theta)-\theta\log\theta}{t+1-\theta}    
% \end{equation}
% attains its global maximum $f_{\theta}(u)\approx 0.434536$ on $[\theta,\infty)$. 
\[\epsilon \coloneqq (\log\log{x})^{-1/2}\quad\text{and}\quad\theta\coloneqq0.4736,\]
and define $u\approx1.2694$ to be the unique point at which the function
\begin{equation}\label{eq:f_{theta}}
f_{\theta}(t)\colonequals\frac{1}{t+1-\theta}\left(\theta\log\frac{t}{\theta}-\left(t-\theta\right)\log\left(1-\frac{\theta}{t}\right)\right)=\frac{t\log t-(t-\theta)\log(t-\theta)-\theta\log\theta}{t+1-\theta}    
\end{equation}
attains its global maximum $f_{\theta}(u)\approx0.4669$ on $[\theta,\infty)$. 

Below, we describe how Proposition \ref{prop:harman} can be used to find a positive integer $k\mid \prod_{p \le (u-\epsilon)\log{x}} p$, along with a set $\mathcal{D}$ of divisors of $k$, where each $d \in \mathcal{D}$ has the property that
\begin{equation}\label{eq:pilowerharman} \pi(x;d,1) \gg \frac{x}{\phi(d)\log{x}}. \end{equation}
Furthermore, $\mathcal{D}$ will be selected in such a way that each $d \in \mathcal{D}$ obeys the estimate
\begin{equation}\label{eq:dsizesecond} d= x^{\theta-\epsilon} \exp\left(O\left(\frac{\log{x}}{(\log\log{x})^2}\right)\right),\end{equation}
and such that
\begin{equation}\label{eq:Dsizesecond}\#\Dd \ge \exp\left(\left(C'+o(1)\right)\frac{\log{x}}{\log\log{x}}\right)=\exp\left(\left(f_{\theta}(u)+o(1)\right)\frac{\log{Y}}{\log\log{Y}}\right),\end{equation}
with \[ Y\coloneqq x^{u+1-\theta}, \quad\text{and}\quad  C'\coloneqq \theta\log\frac{u}{\theta}-\left(u-\theta\right)\log\left(1-\frac{\theta}{u}\right)=(u+1-\theta)f_{\theta}(u).\]

After $k$ and $\mathcal{D}$ have been located, the rest of the argument can be carried out as before. For each $d \in \mathcal{D}$, let $A_d$ denote the number of pairs $(m,p)$ with $m\le y_d=x^{u}/d$ and $p\le x$, satisfying the conditions in \eqref{eq:(m,p)}. Then we have 
\[A_d\gg\pi(x;d,1)\frac{\varphi(d)}{d}\cdot\frac{x^u}{k} \gg \frac{x^{u+1}}{kd\log{x}} = \frac{x^{u+1-\theta+\epsilon}}{k\log{x}}  \exp\left(O\left(\frac{\log{x}}{(\log\log{x})^2}\right)\right),\]
and 
\[ \sum_{d\in\Dd}A_d\gg\frac{x^{u+1-\theta+\epsilon}}{k\log x} \exp\left(\left(f_{\theta}(u)+o(1)\right)\frac{\log{Y}}{\log\log{Y}}\right).  \]
Each pair $(m,p)$ counted by some $A_d$ corresponds to a multiple $m(p-1)$ of $k$ for which \[ m(p-1)\le x^{u+1-\theta+\epsilon} \exp\left(O\left(\frac{\log{x}}{(\log\log{x})^2}\right)\right).\] Comparing the last two displays, we conclude that there is an $n\le x^{u+1-\theta+o(1)}$ with 
\[ 
\omega^{*}(n)\ge \exp\left(\left(f_{\theta}(u)+o(1)\right)\frac{\log{Y}}{\log\log{Y}}\right) \ge\exp\left(\left(f_{\theta}(u)+o(1)\right)\frac{\log{n}}{\log\log{n}}\right).
\]
As $f_{\theta}(u)/\log 2 = 0.67365\ldots > 0.6736$, the estimate \eqref{eq:maxomega*} follows.

To produce $k$ and $\mathcal{D}$, we borrow ideas and results from \cite[pp.\ 647--648]{harman05}. Let 
\[ W \coloneqq \left(\frac25 \log{x}\right)^{3/4}.\] As on p.\ 647 of \cite{harman05}, for some absolute constant $\eta > 0$, there is at most one primitive character $\chi_1\bmod{f}_1$ of conductor \[ f_1 < V\coloneqq \exp(\eta(\log{x})^{3/4}) \] for which $L(s,\chi_1)$ has a zero $\rho$ with   
\begin{equation}\label{eq:zeroregion} \mathrm{Re}(\rho) > 1-\frac{1}{W}, \quad |\mathrm{Im}(\rho)| \le V.\end{equation}
(This follows from the results on exceptional zeros appearing on pp.\ 93--95 of \cite{davenport80}.) We will apply Proposition \ref{prop:harman} with this $\eta$.  Note that if $x^{0.4} < d < x^{\theta}$, in order for \eqref{eq:harmanregion} to fail, $L(s,\chi)$ must have a zero $\rho$ belonging to the region \eqref{eq:zeroregion}. 

If the primitive character $\chi_1\bmod{f_1}$ of the last paragraph exists, we let $p_1$ be a prime factor of $f_1$. Otherwise, we let $p_1=1$. Let $L=(u-\epsilon)\log{x}$ and $R=\pi(L)$ (as before), and set 
\[\rho=\frac{\theta-\epsilon}{u-\epsilon}. \] 
We take
\[k=\prod_{\substack{p\le L\\p\ne p_1}}p=x^{u-\epsilon}\exp\left(O\left(\frac{\log{x}}{\log\log{x}}\right)\right),\]
and we let $d = \prod_{r \mid k} r^{v_r}$, where the $v_r$ are i.i.d. Bernoulli random variables with each $\mathbb{P}(v_r=1)=\rho$. 

By our earlier arguments, a random $d$ satisfies both
\begin{equation}\label{eq:dboundgeneral} |\log{d}-(\theta-\epsilon)\log{x}| < L^{2/3} \end{equation}
and
\begin{equation}\label{eq:rboundgeneral} |\Omega(d) - \rho R| < R^{2/3} \end{equation}
with probability $1+o(1)$. (The possibly missing prime $p_1$ has a negligible effect.) Let $\mathcal{D}_0$ be the set of divisors $d$ of $k$ satisfying \eqref{eq:dboundgeneral} and \eqref{eq:rboundgeneral}. We proceed to remove from $\mathcal{D}_0$ those $d$ for which there is a primitive character $\chi\bmod{f}$, $f\mid d$, where \eqref{eq:harmanregion} fails. Since $\gcd(d,p_1)=1$, in these cases we necessarily have $f \ge V$. Furthermore, $L(s,\chi)$ has a zero in the region \eqref{eq:zeroregion}.

Zero density estimates (e.g., \cite[Theorem 1]{montgomery69} suffices here) show that there are no more than $\exp(O((\log{x})^{1/4}))$ primitive characters $\chi\bmod{f}$, $f\le x$, having a zero in the region \eqref{eq:zeroregion}. (Compare with pp.\ 647--648 of \cite{harman05}.) From that set of characters, throw away those of conductors smaller than $V$, and collect their remaining conductors in a set $\mathcal{F}$. Then each $d$ to be removed from $\mathcal{D}_0$ is divisible by some $f\in \mathcal{F}$. 

We fix $f\in \mathcal{F}$ and examine the probability that $f\mid d$. If $f\nmid k$, then $\mathbb{P}(f\mid d)=0$. Otherwise, $\mathbb{P}(f\mid d) = \rho^{\omega(f)}$, where $\omega(f)$ denotes the number of distinct prime factors of $f$. Since
\[ V \le f \le (u\log{x})^{\omega(f)}, \] we have
\[ \omega(f) \ge \frac{\log{V}}{\log(u\log{x})} = \frac{\eta (\log{x})^{3/4}}{\log(u\log{x})}. \]
Thus (for large $x$),
\[ \mathbb{P}(f\mid d) = \rho^{\omega(f)} < \exp(-(\log{x})^{7/10}). \]
Hence, 
\[ \mathbb{P}(f\mid d\text{ for some $f\in \mathcal{F}$}) \le \#\mathcal{F}\exp(-(\log{x})^{7/10}) = o(1), \]
recalling for the last equality that $\#\mathcal{F} \le \exp(O((\log{x})^{1/4}))$. 

Therefore, after removing all $d\in \mathcal{D}_0$ divisible by an $f \in \mathcal{F}$, we are left with a set $\mathcal{D}$ of divisors of $k$ for which $\mathbb{P}(d\in \mathcal{D}) = 1+o(1)$. We see from \eqref{eq:dboundgeneral} that each $d \in \mathcal{D}$ satisfies \eqref{eq:dsizesecond}. Furthermore, invoking \eqref{eq:rboundgeneral} and repeating the argument leading to \eqref{eq:sizedlower}, we arrive at \eqref{eq:Dsizesecond}. Finally, Proposition \ref{prop:harman} furnishes the desired lower bound \eqref{eq:pilowerharman} on $\pi(x;d,1)$ for all $d \in \mathcal{D}$. 

\section{Concluding remarks}
We have seen that a key, common ingredient in Prachar's argument and the proof of Theorem \ref{thm:maxomega*} is an inequality of the form
\begin{equation}\label{eq:pi(x;d,1)}
\pi(x;d,1)\gg\frac{x}{\phi(d)\log x}  
\end{equation}
for $d\mid k$, where $k$ is essentially the product of all primes $p\le \theta\log x$ with some $\theta\in(0,1)$. The proof of Theorem \ref{thm:maxomega*} reveals that if \eqref{eq:pi(x;d,1)} holds for some fixed $\theta\in(0,1)$, then we have
\[\omega^{*}(n)\ge\exp\left(\left(\max_{t\ge\theta}f_{\theta}(t)+o(1)\right)\frac{\log{n}}{\log\log{n}}\right)\]
for infinitely many $n$, where $f_{\theta}(t)$ is defined as in \eqref{eq:f_{theta}}. In particular, the conjecture of Adleman, Pomerance and Rumely \cite{APR83} mentioned in the introduction, that
\begin{equation}\label{eq:APRconj}
\omega^{*}(n)\ge\exp\left(\left(\log 2+o(1)\right)\frac{\log{n}}{\log\log{n}}\right)    
\end{equation}
for infinitely many $n$, would follow if \eqref{eq:pi(x;d,1)} holds for any fixed $\theta\in(0,1)$ (since, for instance, $\max\limits_{t\ge \theta} f_{\theta}(t) \ge f_{\theta}(2\theta) = \frac{2\theta}{\theta+1}\log{2}$). As \cite[Theorem 2.1]{AGP94} shows, the lower bound \eqref{eq:pi(x;d,1)} is intimately related to zero densities for Dirichlet $L$-functions. An implication of this is that our inequality \eqref{eq:GRHmaxomega*} still holds under the Density Hypothesis which is weaker than GRH. 

Note that in Prachar's argument, if we sum $A_d$ over all divisors $d$ of $k$ instead, we would have
\[\sum_{d\mid k}A_d\gg\frac{x}{k}\sum_{d\mid k}\phi(d)\pi(x;d,1)=\frac{x}{k}\sum_{p\le x}\sum_{\substack{d\mid k\\d\mid p-1}}\phi(d)=\frac{x}{k}\sum_{p\le x}\gcd(p-1,k).\]
The last sum hints at the connection between the maximal order of $\omega^*$ and the distribution of smooth shifted primes $p-1$. Given $y\ge1$, we say that $n\in\N$ is {\it $y$-smooth} if $P^+(n)\le y$, where $P^+(n)$ denotes the greatest prime factor of $n$, with the convention that $P^+(1)=1$. In other words, $y$-smooth numbers are precisely those integers with no prime factors exceeding $y$. For $x\ge y\ge1$, we define the counting functions
\begin{align*}
\Psi(x,y)&\coloneqq\#\left\{n\le x\colon P^+(n)\le y\right\},\\
\pi(x,y)&\coloneqq\#\left\{p\le x\colon P^+(p-1)\le y\right\}.
\end{align*}
In contrast to $\Psi(x,y)$ whose asymptotic behavior is rather well-understood (see for instance \cite[Chapter III.5]{Ten15}), the function $\pi(x,y)$ has remained elusive. Nevertheless, it is widely believed that smooth shifted primes have the same asymptotic density relative to shifted primes as smooth integers do relative to integers. Indeed, Pomerance \cite{Pom80} has conjectured that if $x\ge y\ge1$, then
\begin{equation}\label{eq:smooth p-1}
\frac{\pi(x,y)}{\pi(x)}\sim\frac{\Psi(x,y)}{x}   
\end{equation}
as $y\to\infty$. We conclude our paper with a demonstration that the Adleman--Pomerance--Rumely conjecture \eqref{eq:APRconj} is an easy consequence of Pomerance's conjecture \eqref{eq:smooth p-1}.

Assume \eqref{eq:smooth p-1}. Fix $v>0$ and set $y=v\log x$ with $x$ sufficiently large. The number of pairs $(m,p)$, with $m,p\le x$, $P^+(m)\le y$, and $P^+(p-1)\le y$, is $\Psi(x,y)\pi(x,y)$. Since each $n=(m-1)p\le x^2$ is $y$-smooth, and since the number of $y$-smooth numbers up to $x^2$ is precisely given by $\Psi(x^2,y)$, we deduce from \eqref{eq:smooth p-1} that there is some large $y$-smooth number $n\le x^2$ with at least
\begin{equation}\label{eq:psipi/psi}
\frac{\Psi(x,y)\pi(x,y)}{\Psi(x^2,y)}\sim\frac{\Psi(x,y)^2}{\Psi(x^2,y)}\cdot\frac{\pi(x)}{x}\sim\frac{\Psi(x,y)^2}{\Psi(x^2,y)}\cdot\frac{1}{\log x}   
\end{equation}
representations as $n=m(p-1)$. By \cite[Theorem III.5.2]{Ten15}, we have
\begin{align*}
\log\Psi(x,y)&=\left(1+O\left(\frac{1}{\log\log x}\right)\right)\frac{\log x}{\log y}\int_{0}^{1}\log\left(1+\frac{y}{t\log x}\right)\,\text{dt}\\
&=\left(1+O\left(\frac{1}{\log\log x}\right)\right)\frac{\log x}{\log\log x}\int_{0}^{1}\log\left(1+\frac{v}{t}\right)\,\text{dt},    
\end{align*}
and analogously,
\[\log\Psi(x^2,y)=\left(1+O\left(\frac{1}{\log\log x}\right)\right)\frac{2\log x}{\log\log x}\int_{0}^{1}\log\left(1+\frac{v}{2t}\right)\,\text{dt}.\]
Since $\log(1+z)\ge z/(1+z)$ for all $z>-1$, it follows that
\begin{align*}
\log\frac{\Psi(x,y)^2}{\Psi(x^2,y)}&=\left(\log 2+\int_{0}^{1}\log\left(1-\frac{t}{2t+v}\right)\,\text{dt}+O\left(\frac{1}{\log\log x}\right)\right)\frac{2\log x}{\log\log x}\\
&\ge\left(\log 2-\int_{0}^{1}\frac{t}{t+v}\,\text{dt}+O\left(\frac{1}{\log\log x}\right)\right)\frac{2\log x}{\log\log x}\\
&\ge\left(\log 2-\frac{1}{1+v}+O\left(\frac{1}{\log\log x}\right)\right)\frac{2\log x}{\log\log x}.
\end{align*}
Inserting this in \eqref{eq:psipi/psi}, we find that this particular $n$ has at least
\begin{align*}
\frac{\Psi(x,y)\pi(x,y)}{\Psi(x^2,y)}&\gg\exp\left(\left(\log 2-\frac{1}{1+v}+O\left(\frac{1}{\log\log x}\right)\right)\frac{2\log x}{\log\log x}\right)\\
&\ge\exp\left(\left(\log 2-\frac{1}{v}\right)\frac{\log n}{\log\log n}\right)
\end{align*}
representations as $n=m(p-1)$. Since $v>0$ is arbitrary, this verifies our claim that \eqref{eq:APRconj} follows from \eqref{eq:smooth p-1}.

\bibliographystyle{amsplain}
\bibliography{SPD}

\begin{bibdiv}
\begin{biblist}

\bib{APR83}{article}{
      author={Adleman, {L.\,M.}},
      author={Pomerance, C.},
      author={Rumely, {R.\,S.}},
       title={On distinguishing prime numbers from composite numbers},
        date={1983},
     journal={Ann. of Math.},
      volume={117},
       pages={173\ndash 206},
}

\bib{AGP94}{article}{
      author={Alford, {W.\,R.}},
      author={Granville, A.},
      author={Pomerance, C.},
       title={There are infinitely many {C}armichael numbers},
        date={1994},
     journal={Ann. of Math.},
      volume={139},
       pages={703\ndash 722},
}

\bib{davenport80}{book}{
      author={Davenport, H.},
       title={Multiplicative number theory},
     edition={2nd ed.},
      series={Graduate Texts in Mathematics},
   publisher={Springer-Verlag, New York-Berlin},
        date={1980},
      volume={74},
        note={Revised by Hugh L. Montgomery},
}

\bib{Ding23}{article}{
      author={Ding, Y.},
       title={On a conjecture of {M.\,R.} {M}urty and {V.\,K.} {M}urty},
        date={2023},
     journal={Canad. Math. Bull.},
      volume={66},
       pages={679\ndash 681},
}

\bib{EPS91}{article}{
      author={Erd\H{o}s, P.},
      author={Pomerance, C.},
      author={Schmutz, E.},
       title={Carmichael's lambda function},
        date={1991},
     journal={Acta Arith.},
      volume={58},
       pages={363\ndash 385},
}

\bib{Fan25}{article}{
      author={Fan, {K.\,(S.)}},
       title={The shifted prime-divisor function over shifted primes},
        date={2025},
     journal={Ramanujan J.},
      volume={66},
       pages={1\ndash 46},
}

\bib{FP24}{article}{
      author={Fan, {K.\,(S.)}},
      author={Pomerance, C.},
       title={Shifted-prime divisors},
        date={2024},
     journal={preprint},
        note={arXiv:2401.10427. To appear in a Springer volume dedicated to Helmut Maier on his 70th birthday and edited by J. Friedlander, C. Pomerance, and M. Rassias.},
}

\bib{Gab25}{article}{
      author={Gabdullin, {M.\,R.}},
       title={Moments of the shifted prime divisor function},
        date={2025},
     journal={preprint},
        note={arXiv:2505.24050},
}

\bib{HW08}{book}{
      author={Hardy, {G.\,H.}},
      author={Wright, {E.\,M.}},
       title={An introduction to the theory of numbers},
     edition={6},
   publisher={Oxford University Press, Oxford},
        date={2008},
        note={Revised by D.\,R. Heath-Brown and J.\,H. Silverman},
}

\bib{harman05}{article}{
      author={Harman, G.},
       title={On the number of {C}armichael numbers up to {$x$}},
        date={2005},
     journal={Bull. London Math. Soc.},
      volume={37},
       pages={641\ndash 650},
}

\bib{harman08}{article}{
      author={Harman, G.},
       title={Watt's mean value theorem and {C}armichael numbers},
        date={2008},
     journal={Int. J. Number Theory},
      volume={4},
       pages={241\ndash 248},
}

\bib{montgomery69}{article}{
      author={Montgomery, {H.\,L.}},
       title={Zeros of {$L$}-functions},
        date={1969},
     journal={Invent. Math.},
      volume={8},
       pages={346\ndash 354},
}

\bib{MV06}{book}{
      author={Montgomery, {H.\,L.}},
      author={Vaughan, {R.\,C.}},
       title={Multiplicative number theory. {I. Classical theory}},
      series={Cambridge Studies in Advanced Mathematics},
   publisher={Cambridge University Press, Cambridge},
        date={2006},
      volume={97},
}

\bib{MM21}{article}{
      author={Murty, {M.\,R.}},
      author={Murty, {V.\,K.}},
       title={A variant of the {H}ardy--{R}amanujan theorem},
        date={2021},
     journal={Hardy--Ramanujan J.},
      volume={44},
       pages={32\ndash 40},
}

\bib{Pom80}{article}{
      author={Pomerance, C.},
       title={Popular values of {E}uler's function},
        date={1980},
     journal={Mathematika},
      volume={27},
       pages={84\ndash 89},
}

\bib{Prachar55}{article}{
      author={Prachar, K.},
       title={{\"U}ber die {A}nzahl der {T}eiler einer nat\"urlichen {Z}ahl, welche die form $p-1$ haben},
        date={1955},
     journal={Monatsh. Math.},
      volume={59},
       pages={91\ndash 97},
}

\bib{Ten15}{book}{
      author={Tenenbaum, G.},
       title={Introduction to analytic and probabilistic number theory},
     edition={3},
      series={Graduate Studies in Mathematics},
   publisher={American Mathematical Society, Providence, RI},
        date={2015},
      volume={163},
}

\end{biblist}
\end{bibdiv}
\end{document}